\documentclass[12pt]{amsart}

\usepackage{amsmath,amssymb,amscd,amsthm,amsxtra}
\usepackage[dvips]{graphics,epsfig}

\newtheorem{theorem}{Theorem}
\newtheorem{lemma}{Lemma}
\newtheorem{proposition}{Proposition}

\newcommand{\q}{\quad}

\newcommand{\qqq}{\quad\quad\quad}

\newcommand{\nf}{\infty}

\newcommand{\ga}{\gamma}

\newcommand{\la}{\lambda}

\newcommand{\om}{\omega}
\newcommand{\Om}{\Omega}

\newcommand{\rn}{{\mathbf R}^n}

\newcommand{\lp}{L^{p}}

\newcommand{\li}{L^{\infty}}

\newcommand{\cm}{\mathcal M}

\newcommand{\lab}{\label}

\newcommand{\intl}{\int\limits}

\newcommand{\f}{\frac}
\newcommand{\df}{\displaystyle\frac}

\newcommand{\wt}{\widetilde}
\newcommand{\wh}{\widehat}

\newcommand{\rrr}{\mathbf R}

\begin{document}

\title
{$L^p$ bounds for a maximal dyadic sum operator}

\author{Loukas Grafakos}

\address{
Loukas Grafakos\\
Department of Mathematics\\
University of Missouri\\
Columbia, MO 65211, USA}

\email{loukas@math.missouri.edu}

\author{Terence Tao}

\address{
Terence Tao\\
Department of Mathematics\\
University of California, Los Angeles\\
Los Angeles, CA 90024, USA}

\email{tao@math.ucla.edu}

\author{Erin Terwilleger }

\address{Erin Terwilleger\\
School of Mathematics\\
Georgia Institute of Technology\\
Atlanta, GA 30332 USA}

\email{erin@math.gatech.edu }

\thanks{Grafakos is supported by the NSF.  Tao is a Clay Prize Fellow and is 
supported by a grant from the Packard Foundation.}

\date{\today}

\subjclass{Primary  42A20. Secondary  42A24}

\keywords{Fourier series, almost every convergence}

\begin{abstract}
The authors prove $L^p$ bounds in the range $1<p<\infty$ for a maximal
dyadic sum operator on $\rn$. 
This maximal operator provides a discrete multidimensional model
of   Carleson's operator. Its boundedness is obtained by a simple twist of the 
proof of Carleson's theorem given by Lacey and Thiele \cite{LT} adapted in  
higher dimensions \cite{em}. 
In dimension one,  the $\lp$ boundedness 
of this maximal dyadic sum implies   in 
particular   an alternative proof of 
Hunt's   extension \cite{hunt} of Carleson's theorem  on almost 
everywhere convergence of Fourier integrals. 
\end{abstract}

\maketitle

\section{The Carleson-Hunt theorem}\label{s1}

A celebrated theorem of Carleson \cite{carleson}
states that the  Fourier series of a square-integrable 
function on the circle converges almost everywhere to 
the function. Hunt \cite{hunt} extended this theorem to 
$\lp$ functions for $1<p<\nf$. Alternative proofs of 
Carleson's theorem  were provided by C. Fefferman \cite{fefferman}  
and   by Lacey and Thiele \cite{LT}. The last authors 
proved the  theorem on the line, i.e. they 
showed that for $f$ in $L^2(\rrr)$ the sequence of functions  
$$
S_N(f)(x)=\int_{|\xi |\le N} \wh{f}( \xi) e^{2\pi i x\xi} d\xi 
$$
converges to $f(x)$ for almost all $x\in \rrr$ as $N\to \nf$. 
This result was obtained as a consequence of the 
boundedness of    the maximal operator  
$$
\mathcal C(f) = \sup_{N>0} |S_N(f)|
$$
from $L^2(\rrr)$ into $L^{2,\nf}(\rrr)$. 
In view of the transference theorem 
of Kenig and Tomas \cite{kenig-tomas} the above result is 
equivalent to the analogous theorem for Fourier series on the circle. 
Lacey and Thiele \cite{LT-p} have also obtained a proof of Hunt's 
theorem by adapting the techniques in  \cite{LT} to the $L^p$ case
but this proof is rather complicated compared with the 
relatively  short and elegant proof they gave for $p=2$.

Investigating higher dimensional analogues, 
Pramanik and Terwilleger \cite{em} recently adapted the proof of
Carleson's theorem by Lacey and Thiele
\cite{LT} to prove weak type $(2,2)$ bounds for a discrete maximal operator on
$\rn$ similar to the one which arises in the aforementioned proof.
After a certain averaging procedure, this
result provides an alternative proof of Sj\"olin's \cite{sjolin}
theorem on the weak $L^2$ boundedness of
  maximally modulated Calder\'on-Zygmund operators on $\rn$.
The purpose of this note is to extend the result of Pramanik and Terwilleger
\cite{em} to the range $1<p<\infty$ via a   variation of the 
$L^2\to L^{2,\nf}$ case. 
Particularly in   dimension $1$, the theorem below
yields a new   proof of 
Hunt's theorem (i.e. the $\lp$ boundedness of $\mathcal C$   
for $1<p<\nf$) using a variation of the proof of Lacey and Thiele
\cite{LT}.

\section{Reduction to two estimates}\label{s2}

We use the notation introduced in 
\cite{LT} and expanded in \cite{em}. A tile in $\rn \times \rn$ is a  product 
of dyadic cubes of the form
$$
 \prod_{j=1}^n I^j = \prod_{j=1}^n \, \, [m_j 2^k,(m_j +1) 2^k),
 $$
 where  $k$ and $m_j$ are integers for all $j= 1,2, \cdots, n$. 
We denote a {\em tile}  by $s=I_s\times \om_s$, where $|I_s||\om_s|=1$. 
The cube $I_s$ will be called the
time projection of
$s$ and  $\om_s$ the frequency projection of $s$.
For a  tile $s$ with $\om_s=\om_s^1 \times \om_s^2 \times \ldots
\times \om_s^n $, we can divide each dyadic interval $\om_s^j$ into two
intervals of the form
$$
\om_s^j= (\om_s^j \cap (-\nf,c(\om_s^j)) \cup 
(\om_s^j \cap [c(\om_s^j),\nf))
$$
for $j=1,2, \ldots, n$.
 Then $\om_s$ can be decomposed into 
$2^n$ subcubes formed from all combinations
of cross products of these half intervals.  We number these subcubes using the
lexicographical order on the centers and denote the subcubes by $\om_{s(i)}$
for $i=1,2, \ldots, 2^n$.  A tile $s$ is then the union of $2^n$
{\em {semi-tiles}} given by  
$s(i)=I_s \times \om_{s(i)}$ for $i=1,2, \ldots, 2^n$. 

We let $\phi$ be a Schwartz function  such that
$\wh{\phi}$ is real, nonnegative, and supported in the cube 
$[-1/10, 1/10]^n$. Define
\begin{equation*} 
\phi_{s} (x)=  |I_s|^{-\f12} 
\phi\bigg(\f{x-c(I_s)}{|I_s|^{\f{1}{n}}} 
\bigg) e^{2\pi i c(\om_{s(1)}) \cdot x}\, ,   
\end{equation*}
where 
$c(J)$ is the center of a cube $J$.
As in \cite{LT} and \cite{em}, we will consider 
the dyadic sum operator 
$$
\mathcal D_r(f)= \sum_{s\in D} \langle f, \phi_s\rangle (\chi_{\om_{s(r)}}\circ
N)\phi_s\, ,
$$
where $2 \le r \le 2^n$ is a fixed integer, 
$N:\,\rn\to \rn$ is a fixed measurable function, $D$ is a 
set of tiles, and  $\langle f, g\rangle$ is the complex inner product 
$\int_{\rrr} f(x)\overline{g(x)}\, dx$.

The following theorem is the main result of this article. 
\begin{theorem}
Let $1<p<\nf$. Then there is a constant   $C_{n,p}$   independent of the 
measurable function $N$, of the set $D$, and of $r$ such that 
for all $f\in \lp(\rn)$ we have 
\begin{equation}\lab{2}
\| \mathcal D_r(f) \|_{\lp(\rn)}\le C_{n,p} \|f\|_{\lp(\rn)} \, . 
\end{equation}
\end{theorem}
In one-dimension, using the averaging procedure introduced  in \cite{LT},
it follows that the norm estimate (\ref{2}) implies 
\begin{equation*}
\| \mathcal C(f) \|_{\lp}\le C_p \|f\|_{\lp}, 
\end{equation*}
which is the Carleson-Hunt theorem.
Using   the Marcinkiewicz interpolation theorem \cite{S}
and the restricted weak type reduction of Stein and Weiss \cite{SW},  
estimate (\ref{2}) will be a consequence of 
the restricted weak type estimate 
\begin{equation}\lab{3}
\| \mathcal D_r(\chi_F) \|_{L^{p,\nf}(\rn)}\le C_{n,p}  |F|^{\f 1p}\, , 
  \qqq 1<p<\nf    
\end{equation}
which is supposed to hold 
for all $n$-dimensional sets $F$ of finite measure. 
But to show that a function $g$ lies in $L^{p,\nf}$, it suffices to 
show that for every measurable set $E$ of finite measure, there is 
a subset $E'$ of $E$ which satisfies $|E'|\ge \f12 |E|$   and also  
$$
\bigg| \int_{E'} g(x) \, dx\bigg|\le   A\, |E|^{\f{p-1}{p }}\, ;
$$
this implies that $\|g\|_{L^{p,\nf}(\rn)}\le c_{ p}\, A$, where $c_{ p}$ is a
constant that depends only on $p$.

Let $C (n,q)$ be the weak type $(q,q)$ operator norm
 for the Hardy-Littlewood maximal operator. 
Given a set $E$ of finite measure  we set 
$$
\Om=\Big\{M(\chi_F)>\Big(2\frac{|F|}{|E|}\Big)^{\!\f1q}C (n,q)\Big\}, 
$$
where we choose $q$ so that  $p < q \le \infty$ if $ |F|  > |E|$ and 
$1 \le q <p$ if $ |F| \le |E|$. Note that in the first case the 
set $\Om$ is empty. 
Using the $L^q$ to $L^{q,\infty}$ boundedness of the Hardy-Littlewood maximal operator,
we have $ |\Om| \le \f{1}{2}|E|$
and hence $|E'|\ge \f12 |E|$. Thus estimate (\ref{3}) will follow from  
\begin{equation}\lab{4}
\bigg|\int_{E'} \mathcal D_r(\chi_F)(x) \, dx \bigg|\le C_{n,p} |E|^{\f{p-1}{p
}} |F|^{\f 1p}\,  ,  
\end{equation}
where $C_{n,p}$ depends only on $p$ and dimension $n$. 
The required estimate (\ref{4}) will then be a consequence of the 
following two  estimates:
\begin{equation}\lab{5}
\bigg|\int_{E'} \sum_{\substack{ s\in D\\ I_s\subseteq \Om}} \langle \chi_F ,
\phi_s\rangle (\chi_{\om_{s(r)}}\circ N)\phi_s(x) \, dx \bigg|\le 
C_{n,p,q}
|E|^{\f{p-1}{p }} |F|^{\f 1p}\,  ,  
\end{equation}
and 
\begin{equation}\lab{6}
\sum_{\substack{ s\in D\\ I_s\nsubseteq \Om}} 
|\langle \chi_F,\phi_s\rangle |\, |\langle     \chi_{E'\cap N^{-1}[\om_{s(r)}]} ,
\phi_s \rangle |\le 
 C_{n,p,q}
|E|^{\f{p-1}{p }} |F|^{\f 1p}\,  .  
\end{equation}

\section{The proof of estimate (\ref{5})}

Following \cite{LT2}, we denote by $I(D)$   the dyadic grid  
which consists of all the
time projections of tiles in $D$. For each dyadic cube $J$ in $I(D)$  we define
$$
D_J:= \{s \in D:I_s=J\}
$$
and a function 
$$
\psi_J (x):=|J|^{-\f12} \bigg(1+\f{|x-c(J)|}{|J|^{\f{1}{n}}}\bigg)^{-\ga},
$$
where $\ga$ is a large integer to be chosen shortly.  
For each $k=0,1,2,\dots$ we introduce families 
$$
\mathcal F_k  =\big\{ J \in  I(D) :\,\,  
2^k J\subseteq \Om ,\,\, 2^{k+1}J\nsubseteq \Om  \big\}.
$$
We may assume $|F|\le |E|$, otherwise the set $\Om$ is empty and  (\ref{5}) is
trivial.

We begin by controlling
the left hand side of (\ref{5})  by
\begin{align}\begin{split}\lab{10.6.erin1}
 &\,\,\sum_{\substack{ J\in  I(D)\\ J\subseteq \Om }}
 \bigg | \sum_{\substack{ s \in D(J)}}
 \int_{E'} \big\langle \chi_F  \, | \,
\phi_s\big\rangle  \,   \chi_{\om_{s(2)}}( N (x) )\phi_s (x) \,dx \bigg|    \\
\le  &\,\,\sum_{k=0}^\nf \sum_{\substack{J\in   I(D)\\ J\in\mathcal F_k} } 
\bigg|\int_{E'} \sum_{s\in D(J)} 
\big\langle \chi_F \, | \, \phi_s\big\rangle  \, \chi_{\om_{s(r)}}( N (x) )  \phi_s(x)
\, dx  \bigg| 
\end{split}\end{align}

Using the fact that the function $M(\chi_F)^{\f12}$ is an $A_1$ weight with 
$A_1$-constant bounded above by a quantity independent of $F$, it is easy to 
find  a  constant  $C_0<\nf$ 
such that for each $k=0,1,\dots$ 
and $J\in \mathcal F_k$ we   have 
\begin{equation}\lab{10.6.H.7}
 \big\langle \chi_F , \psi_J\big\rangle  
 \le   |J|^{\f12}\inf_{J}
M(\chi_F)   
 \le   |J|^{\f12}\, C_0^k\,\inf_{2^{k+1}J} M(\chi_F)  
 \le  
C(n,q)2^{\f1q}\,C_0^k |J|^{\f12} \Big(\f{|F|}{|E|}\Big)^{\f1q}\, 
\end{equation}
since $2^{k+1}J$ meets the complement of $\Om$. 
For $J\in \mathcal F_k$ one also has that $E'\cap  2^kJ =\emptyset$ and hence 
\begin{equation}\lab{10.6.H.8}
\int_{E'} \psi_J(y) \, dy \le \int_{(2^kJ)^c} \psi_J(y) \, dy
\le |J|^{\f12} C_\ga 2^{-k\ga}\, . 
\end{equation} 

Next we note that for each $J\in I(D)$ and $x\in \rn$  there is at 
most one $s=s_x\in D_J$ such that $N(x)\in \om_{s_x(r)}$. 
Using this observation along with 
(\ref{10.6.H.7}) and (\ref{10.6.H.8}) we  can therefore  
estimate  the  expression  on the right in (\ref{10.6.erin1})   as follows
\begin{align}
\le &\,\,\sum_{k=0}^\nf \sum_{\substack{J\in   I(D)\\ J\in\mathcal F_k} } 
\bigg|\int_{E'}  
\big\langle \chi_F \, | \, \phi_{s_x }\big\rangle  \, \chi_{\om_{s_x(r)}}( N (x) ) 
\phi_{s_x}(x) \, dx  \bigg| \notag \\
\le &\,\,C  \sum_{k=0}^\nf \sum_{\substack{J\in   I(D)\\ J\in\mathcal
F_k} }  \int_{E'}  
\big\langle \chi_F , \psi_{J}\big\rangle   
\psi_{J}(x t) \, dx   \notag \\
\le &\,\,C  \,
\Big(\f{|F|}{|E|}\Big)^{\f1q}\sum_{k=0}^\nf C_0^k   
\sum_{J\in \mathcal F_k } |J|^{\f12} \, \int_{E'} \psi_J(x)\, dx  \notag \\ 
\le &\,\,  C  \,\Big(\f{|F|}{|E|}\Big)^{\f1q}\sum_{k=0}^\nf (C_02^{-\ga})^k   
\sum_{J\in \mathcal F_k } |J| \,  \lab{10.6.erin2}  
\end{align}  
and at this point we pick $\ga$ so that $C_02^{-\ga}<1$. It remains to control
$ 
\sum_{J\in \mathcal F_k } |J| 
$ 
for each nonnegative integer $k$. In doing this we let $\mathcal F_k^*$ 
be all elements of $\mathcal F_k$ which are maximal under inclusion. Then 
we observe that if $J\in \mathcal F_k^*$  and $J'\in \mathcal F_k $ 
satisfy $J'\subseteq J$ then $\text{dist }(J',J^c)=0$  (otherwise 
$2J'$ would be contained in $J$ and 
thus $2^{k+1}J'\subseteq 2^kJ\subseteq \Om $.) But for any fixed 
$J$ in $\mathcal F_k^*$  and any scale $ m$,
all the cubes $J'$ in $J'\in \mathcal F_k $ of sidelength $2^m$ that 
touch $J$ are concentrated near the boundary of $J$ and 
have total measure at most $2^m \cdot 2^n (|J|^{\f1n})^{n-1}$. 
Summing over all integers $m$ with $2^m\le |J|^{\f1n}$, we obtain a bound which is 
at most a multiple of $|J|$.   We conclude that 
$$
\sum_{J\in \mathcal F_k } |J|  
= \sum_{J\in \mathcal F_k^* }   \sum_{\substack{J'\in \mathcal F_k
\\  J'\subseteq J}} |J'| \le 
\sum_{J\in \mathcal F_k^* } c_n \, |J| \le c_n \, |\Om |
$$
since elements of $\mathcal F_k^* $ are disjoint and contained in $\Om $. 
Inserting this estimate in (\ref{10.6.erin2}) 
and using that the Hardy-Littlewood maximal operator is of weak type $(1,1)$,
we obtain the required  bound
$$
 C \Big(\f{|F|}{|E|}\Big)^{\f1q}\, |\Om | \le C'  \, |F|  
 \le   C'  \, |E|^{\f{p-1}{p }} |F|^{\f 1p}  
$$
for the    
expression  on the right in (\ref{10.6.erin1}) 
and hence for the expression on the left in (\ref{5}).

\section{The proof of estimate (\ref{6})}

In proving estimate  (\ref{6}) we may assume that 
$\f12 \le |E|\le 1$ by a simple scaling argument. 
(The scaling changes the sets   $D$, $\Om$, and the measurable function 
$N$ but note that  the final constants are independent of these 
quantities.) In addition all constants in the sequel are allowed to depend
on $n$ and $ p $  as described above.
We may also assume that the set $D$ is finite. 
Note that under the normalization of the set $E$,  our choice of $q$ 
is as follows: $1\le q<p$ if $|F|\le c_0$ and $p<q\le \nf$ when $|F|> c_0$
where $c_0$ is a fixed number in the interval $(\f12, 1)$, (in fact $c_0=|E|$)  

We recall that 
a finite set of tiles $ T$ is called a tree if there 
exists a tile $t\in  T$ such that all $s\in  T$ satisfy $s<t$ 
(which means $I_s\subset I_{t}$ and $ \om_{t}\subset \om_s$.)
In this case we call $t$   the top of $ T$ and we  
denote it by $t= t(  T)$.
A   tree $ T$   is called an $r$-tree if 
$$
\om_{t( T)(r)} \subset \om_{s(r)}
$$ 
for all $s\in  T$. 
For a finite set of tiles $Q$ we define the energy of a nonzero function 
$f$ with respect to $Q$ by 
$$
 \mathcal E( f; Q) = \f{1}{\|f\|_{L^2(\rn)}} \sup_{  T} 
\bigg(\f{1}{|I_{ t(  T)}| } \sum_{s\in   T} 
  |\langle f, \phi_s\rangle |^2 \bigg)^{\f12}\, , 
$$
where the supremum is taken over all $r$-trees $ T$ contained in 
$ Q$.
We also define the mass of a set of tiles $Q$ by 
$$
\cm(Q)= \sup_{s \in Q}\sup_{\substack {u \in Q \\ s<u }} \ 
\intl_{E' \cap N^{-1}[\om_{u(r)}]} \df {|I_u|^{-1}} 
{\left(1+ \f {|x-c(I_u)|} {|I_u|^{1/n}}\right)^{\!\gamma}} \, dx.
$$

We now fix a set of tiles $D$ and sets 
$E$ and $F$ with finite measure (recall $\f12\le |E|\le 1$).  
We define $P$ to be the set of all 
tiles in $D$ with the property $I_s\nsubseteq \Om$. 
Given a finite set of tiles $ P$, 
find   a very large integer $m_0$ one can construct a
sequence of  pairwise disjoint sets  $P_{m_0}$, 
$P_{m_0-1}$, $P_{m_0-2}$,  $P_{m_0-3}$, ... such
that 
\begin{equation*} 
P= \bigcup_{j=-\nf}^{m_0} P_{j}  
\end{equation*}
and such that the following properties are satisfied 
\begin{enumerate}
\item[(a)] $\mathcal E(\chi_F;P_{ j}) \le 2^{(j+1)n} $   for all $j\le m_0$.
\item[(b)] $\mathcal M(P_{ j}) \le 2^{(2j+2)n}$   for all $j\le m_0$.
\item[(c)] $\mathcal E\big(\chi_F;P \setminus (P_{m_0}  
\cup \dots \cup P_{j} )\big)  \le 2^{ jn } $   for all $j\le m_0$.
\item[(d)] $\mathcal M\big(P \setminus (P_{m_0} 
\cup \dots \cup P_{j} )\big)  \le 2^{ 2jn }$   for all $j\le m_0$.
\item[(e)] $P_j$ is a union of trees $T_{jk}$ 
such that 
$\sum_{k} |I_{t ({T_{jk}})}| \le C_0 2^{-2jn}$ 
for all $j\le m_0$.
\end{enumerate}
This can be done by induction, see \cite{fefferman}, \cite{LT}, and 
is based  on an energy and a mass lemma shown in \cite{em}. 

The following lemma is the main ingredient of the proof 
and will be proved in the next section.

\begin{lemma}\lab{L}
There is a constant $C$ such that for all measurable sets $F$ and 
all finite set of tiles $P$ which satisfy 
$I_s\nsubseteq \Om$ for all $s\in P$, 
we have   
$$
\mathcal E(\chi_F; P) \le C |F|^{\f1q-\f12}
$$
\end{lemma}

Note that this gives us decay no matter if $|F|$ is large or small due to the choice
of $q$ (the reader is reminded that if $|F|\le c_0$ then $q\in [1,p)$ while
if $|F|\ge c_0$ then $q\in (p,\nf]$.)
We also recall the estimate below from \cite{em}.

\begin{lemma}\lab{L2}  
There is a finite constant $C_1$ such that for all trees  
 $  T$, all $f\in L^2(\rn)$, and 
all measurable sets $E'$ with $|E'|\le 1 $ we have 

\begin{equation}\lab{10.5.basicestimate}
\sum_{s\in   T}  \, \big|
 \langle  f,   \phi_s  \rangle 
\langle     \chi_{E'\cap N^{-1}[\om_{s(r)}]} , \phi_s \rangle  
\big|  \le C_1 \, |I_{t( T)}|\, \mathcal E(f; T)\, 
\mathcal M( T) \|f\|_{L^2(\rn)}. 
\end{equation}
\end{lemma}

Given the sequence of sets  
$P_j$ as above,  we
use (a), (b), (e), the observation that the mass is always bounded by $1$,
and    Lemmata  \ref{L} and \ref{L2} to obtain 
{\allowdisplaybreaks
\begin{align*}
\,\,&\,\,\sum_{s\in P}  \, \big|
 \langle  \chi_F,   \phi_s  \rangle 
\langle     \chi_{E'\cap N^{-1}[\om_{s(r)}]} , \phi_s \rangle   \big| \\
= \,\,&\,\,\sum_{j} \sum_{s\in P_{j }}  \, \big|
 \langle  \chi_F,   \phi_s  \rangle 
\langle     \chi_{E'\cap N^{-1}[\om_{s(r)}]} , \phi_s \rangle   \big| \\
\le \,\,&\,\,\sum_{j}\sum_{k} \sum_{s\in T_{jk}}  \, \big|
 \langle  \chi_F,   \phi_s  \rangle 
\langle     \chi_{E'\cap N^{-1}[\om_{s(r)}]} , \phi_s \rangle   \big| \\
\le \,\,&\,\,C_1\sum_{j}\sum_{k}  
|I_{t(T_{jk})}|\, \mathcal E
(T_{jk})\, \mathcal M(T_{jk}) |F|^{\f12}\\
\le \,\,&\,\,C_1\, |F|^{\f12}\sum_{j}\sum_{k}  
|I_{t(T_{jk})}|\, \min(2^{(j+1)n}, C|F|^{\f1q-\f12})\, \min (1, 2^{(2j+2)n}) \\
\le \,\,&\,\, C' |F|^{\f12}\sum_{j}   2^{-2jn}  \, 
\min(2^{j n},  |F|^{\f1q-\f12}) \min (1,2^{2jn })  \\
\le \,\,&\,\, C'' |F|^{\f1q}\big(1+\big|\log  |F|^{\f12-\f 1q} \big|\big)\\
\le \,\,&\,\, C''' \min(1,|F|) \big(1+\big|\log  |F|  \big|\big)\\
\le \,\,&\,\, C_p |F|^{\f 1p}\, 
\end{align*}
}$\!\!\!$
for all $1<p<\nf$. We observe   
that the choice of $q$ was made to deal with the logarithmic presence 
in the estimate above. Had we taken $q=p$ throughout, we would have 
obtained the sought estimates with the extra factor of  
$ 1+\big|\log  |F|  \big| $. 

Looking at the penultimate inequality above, we  note that we have actually obtained
a stronger estimate than the one claimed in   (\ref{4}). 
Rescaling the set $E$ and taking $q$ to be either $1$ or $\nf$, 
we   have actually proved
that   for every measurable set $E$ 
of finite measure, there is a 
subset $E'$ of $E$ such that for all measurable sets $F$ of finite measure we have 
$$
\Big|\int_{E'} \mathcal D_r(\chi_F)\, dx \Big|\le C\, |E|\, 
\min\Big(1,\f{|F|}{|E|}\Big) \Big(1+\Big|\log  \f{|F|}{|E|}  \Big|\Big)\, .
$$
This will be of use to us in section 6.  

\section{The proof of Lemma \ref{L}}

It remains to prove Lemma \ref{L}. Because of our normalization of 
the set $E$ we may assume that $\Om=  \{M(\chi_F)> c\, |F|^{\f1q}\}$ for some 
$c>0$. 
Fix an $r$-tree $T$ contained in $P$ and 
let $I_t=I_{t(T)}$ be the time projection of its top. 

We write the function $\chi_F$ as  
$\chi _{F\cap 3I_t}+\chi_{F\cap (3I_t)^c}$. 
We begin by observing  that for $s$ in $P$ one has 
$$
|\langle \chi_{F\cap (3I_t)^c} , \phi_s\rangle | \le \f{C_\gamma |I_s|^{ \f12}
\inf\limits_{I_s} M(\chi_F)} {\Big(1+\df{\text{dist}((3I_t)^c,
c(I_s)}{|I_s|^{\f1n}}\Big)^{ \gamma}}
\le  C_\gamma |I_s|^{ \f12}
|F|^{\f1q}\, \bigg(\f{|I_s|}{|I_t|}\bigg)^{\!\f{\gamma}{n}} 
$$
since  $I_s$ meets the complement of $\Om$ for every $s\in P$. 
Square this inequality and sum over all $s$ in $T$ to obtain 
$$
 \sum_{s\in T}|\langle \chi_{F\cap (3I_t)^c} , \phi_s\rangle |^2
\le C \, |I_t|\, |F|^{\f{2}{q}}\, , 
$$
where the last estimate follows 
by placing the $I_s$ 's into groups $\mathcal G_m$ of cardinality at most
$2^{mn}$
so that each element of  $\mathcal G_m$ has size $2^{-mn}|I_t|$. 

We now turn to the corresponding estimate for the function $\chi _{F\cap
3I_t}$. 
At this point it will be convenient
to distinguish   the   case $|F|> c_0$
from the   case  $|F|\le c_0  $. In the   case $|F|> c_0$ the set $\Om $ 
 is empty and  therefore
$$
 \sum_{s\in   T}|\langle \chi_{F\cap  3I_t } , \phi_s\rangle |^2
\le C\,  \| \chi_{F\cap 3I_t} \|_{L^2}^2 \le 
C \, |I_t|\, \le C\,   |I_t|\, |F|^{   \f 2q}\, ,
$$
where the first estimate follows follows from the 
Bessel inequality (\ref{bessel-0})  which holds on any $r$-tree $T$;
the reader may consult \cite{em} or prove it directly.

We therefore concentrate on the case $|F|\le c_0$. 
In proving   Lemma \ref{L} we may assume that there exists a 
point $x_0\in I_t$ such that $M(\chi_F)(x_0)\le c\, |F|^{\f1q}$, otherwise 
there is nothing to prove. 
We may also assume that  the center of 
$\om_{t(T)}$ is zero, i.e. $c(\om_{t(T)})=0$, otherwise we may work with a 
suitable modulation of the function $\chi_{F\cap 3I_t}$ in
the Calder\'on-Zygmund decomposition below.

We write the set $\Om=\{M(\chi_{F})> c\, |F|^{\f1q}\}$ as a disjoint 
union of dyadic cubes $J_\ell'$ such that 
the dyadic parent $\wt{J_\ell'}$ of $J_\ell'$ is not contained in $\Om$ and 
therefore  
$$
|F \cap J_\ell'|\le 
|F \cap  \wt{J_\ell'}|\le 2 \,c \,|F|^{\f1q}\, |J_\ell'|\, . 
$$
Now some of these dyadic cubes may have size larger than or equal to  
 $|I_t|$.  Let $ J_\ell'$ be such a cube. Then we split 
$J_\ell'$ in $\f{|J_\ell '|}{|I_t|}$ cubes $J_{\ell,m}'$ each of size
exactly
$|I_t|$. Since there is  an $x_0\in I_t$ with $M(\chi_F)(x_0)\le c\, |F|^{\f1q}$, it
follows that 
\begin{equation}\lab{11}
|F\cap J_{\ell,m}' |\le 2\,c \,|F|^{\f1q} \, |I_t|\, \bigg(1+\f{\text{dist}
(I_t,J_{\ell,m}')}{|I_t|^{\f1n}}\bigg)^{\!n}\, . 
\end{equation}
We now have a new collection of dyadic cubes $\{J_k\}_k$ contained
in  $\Om$ consisting of all the previous $J_\ell'$ when $|J_\ell'|< |I_t|$
and  the $J_{\ell,m}'$'s when $|J_{\ell,m}'|\ge  |I_t|$. In view of the
construction we have
\begin{equation}\lab{11p}
| F\cap J_k |\le \begin{cases}
2\, c \,|F|^{\f1q} \, |J_k|\,   &\text{when $|J_k|<|I_t|$}\\
2\, c \,|F|^{\f1q} \, |J_k|\, \bigg(1+\df{\text{dist}
(I_t,J_k)}{|I_t|}\bigg)^{\!n}  &\text{when $|J_k|=|I_t|$}
\end{cases}
\end{equation}
for all $k$. 
We now define the ``bad functions''
$$
b_k = \chi_{J_k\cap 3I_t\cap F} - \f{|J_k\cap 3I_t\cap F|}{|J_k|}\chi_{J_k}
$$
which are supported in $J_k$, have mean value zero, and they
satisfy 
$$
\|b_k\|_{L^1(\rn)}\le 2\,c \,|F|^{\f1q} \, |J_k|\, \bigg(1+\f{\text{dist}
(I_t,J_k)}{|I_t|}\bigg)^{\!n}\, . 
$$
We also set 
$$
g= \chi_{F\cap 3I_t} - \sum_{k} b_k
$$
the ``good function'' of the 
above Calder\'on-Zygmund decomposition. We check that 
that $ \|g\|_{L^\nf(\rn)}\le C |F|^{\f1q} $. Indeed, for $x$ in $J_k$ we have 
$$
g(x)= \f{|F \cap   3I_t\cap J_k |}{|J_k|}\le \begin{cases}
\df{|F\cap J_k|}{|J_k|}&\text{when $|J_k|< |I_t|$}\\ \q &\q \\
\df{|F\cap 3I_t|}{|I_t|}&\text{when $|J_k|= |I_t|$}
\end{cases}
$$
and both of the above are at most a multiple of $|F|^{\f1q}$; the latter is 
because there is an $x_0 \in I_t$ with $M(\chi_F)(x_0)\le c \, |F|^{\f1q}$. 
Also for  $x\in (\cup_k J_k)^c=\Om^c$, $g(x)=\chi_{F\cap 3I_t}(x)$ which is 
at most $M(\chi_F)(x)\le c\, |F|^{\f1q}$. We conclude that 
$\|g\|_{\li(\rn)}\le C \, |F|^{\f1q}$.
Moreover 
$$
\|g\|_{L^1(\rn)} \le \sum_k \int_{J_k}\f{|F \cap   3I_t\cap J_k |}{|J_k|}\, dx + 
\|\chi_{F\cap 3I_t}\|_{L^1(\rn)} \le C \, |F\cap 3I_t|\le C\, |F|^{\f1q}\, |I_t|
$$
since the $J_k$ are disjoint. It follows that 
$$
\|g\|_{L^2(\rn)} \le C \, |F|^{\f{1}{2q}} |F|^{\f{1}{2q}}\, 
|I_t|^{\f12} =C |F|^{\f1q}\,
|I_t|^{\f12}\, . 
$$
Using the simple Bessel inequality
\begin{equation}\lab{bessel-0}
\sum_{s\in T} |\langle g, \phi_s\rangle|^2 \le C \, \|g\|_{L^2(\rn)}^2
\end{equation}
we obtain the required conclusion for the function $g$. 

For a fixed $s\in P$ and $J_k$ we will denote by 
$$
d(k,s)= \text{dist } (J_k,I_s)\, .
$$
Then we  have the following  estimate for all $s$ and $k$:
\begin{equation}\lab{122}
|\langle b_k , \phi_s \rangle | \le C_\gamma\, |F|^{\f1q}\, |J_k| 
\Big(1+\f{d(k,t)}{|I_t|^{\f1n}}\Big)^{\!n}\, 
\f{|J_k|\,|I_s|^{-\f32}}{
(1+\f{d(k,s)}{|I_s|^{\f1n}})^{\gamma+n}}  
\le  \f{C_\gamma\, |F|^{\f1q} \,  |J_k|^2\,  |I_s|^{-\f32} }{
(1+\f{d(k,s)}{|I_s|^{\f1n}})^{\gamma}}  
\end{equation}
since $ 1+\f{d(k,t) }{|I_t|^{\f1n}} \le  1
+\f{d(k,s)}{|I_s|^{\f1n}}$. 

We also   have the  estimate
\begin{equation}\lab{142}
|\langle b_k , \phi_s \rangle | \le 
\f{C_\gamma\,   \,|F|^{\f1q}\, |I_s|^{ \f12}}{
(1+\f{d(k,s)}{|I_s|^{\f1n}})^{\gamma}}  \, . 
\end{equation}

To prove   (\ref{122}) we use the fact  that the center of 
$\om_{t(T)}=0$ (which implies that 
$\phi_s'$ obeys   size estimates similar to $|I_s|^{-1}|\phi_s|$) 
and the mean value property of  $b_k$ to obtain 
 $$
\big| \langle b_k , \phi_s \rangle \big|  =\Big|\int_{J_k}
b_k(y)\big( \phi_s(y)-\phi_s(c(J_k))\big)\, dy\Big|\le 
\|b_k\|_{L^1} |J_k| \,  \sup_{\xi\in J_k} 
\f{C_\gamma|I_s|^{-\f32}}{(1+\f{|\xi-c(I_s)|}{|I_s|^{\f1n}})^\gamma }\, . 
$$
To prove estimate (\ref{142}) we note that 
$$
|\langle b_k , \phi_s \rangle | \le C_\gamma\, |I_s|^{ \f12}\, \big(\inf_{I_s}
M(b_k) \big) \f{1}{ (1+\f{d(k,s)}{|I_s|^{\f1n}})^{\gamma}} 
$$
and that 
$$
M(b_k) \le M(\chi_F) + \f{|F \cap 3I_t \cap J_k|}{|J_k|}M(\chi_{J_k})
$$
and since $I_s\nsubseteq \Om$ we have $\inf_{I_s}M(\chi_F) \le c \, |F|^{\f1q}$
while the second term in the sum above was observed earlier to be at 
most $C\, |F|^{\f1q}$. 

Finally we have the estimate
\begin{equation}\lab{132}
|\langle b_k , \phi_s \rangle | 
\le  \f{C_\gamma\, |F|^{\f1q} \,  |J_k| \,  |I_s|^{-\f12} }{
(1+\f{d(k,s)}{|I_s|^{\f1n}})^{\gamma}}  
\end{equation}
which follows by taking the
geometric mean of (\ref{122}) and (\ref{142}).

Now for a fixed $s\in P$ we may have either $J_k\subseteq I_s$ or 
$J_k\cap I_s = \emptyset$ (since $I_s$ is not contained in $\Om$.)
Therefore for   fixed $s\in P$ there are only three possibilities for
$J_k$:\\ (a) $J_k \subseteq 3I_s$\\
(b)  $J_k \cap 3I_s=\emptyset$\\
(c)  $J_k \cap I_s=\emptyset$, $J_k \cap 3I_s\neq \emptyset$, and $J_k
\nsubseteq 3I_s$. \\
Observe that case (c) is equivalent to the following statement: \\
(c) $J_k \cap I_s=\emptyset$, 
$d(k,s)=0$,  and $|J_k|\ge 2^n|I_s|$.

Let us start with case (c).  Note that for each $I_s$ there exists at most 
$2^n-1$ choices
of $J_k$ with the above properties.  Thus  for each $s$ in the sum below 
we can pick one $J_{k(s)}$ at a cost of   $2^n-1$, which is harmless.
Also note that since $d(k,s)=0$  and $|J_k|\ge 2^n|I_s|$, we must have that $I_s
\subset 2J_k$.  But $I_s \subset I_t$ and $|J_k| \le |I_t|$ implies that $J_k
\subset 3I_t$.
Now for a given $J_k$ and a fixed
scale $m \ge 1$, there are at most 
$2^m\times (\# \textup{ of sides}) + 2^n$ possibilities of $I_s$ 
such that $2^{-mn}|J_k|= |I_s|$ and $d(k,s)=0$. 
Using (\ref{142}) we obtain 
\begin{align*}
\sum_{s\in T} \Big|\sum_{k:\,\, J_k\text{ as in (c)}}
\langle b_k, \phi_s\rangle \Big|^2
\le &\,\, (2^n-1)^2 \sum_{s\in T} \Big| 
\langle b_{k_(s)}, \phi_s\rangle \Big|^2 \\
\le &\,\, C_n\, |F|^{\f{2}{q}} \! \sum_{\substack{s\in T\text{ for which}\\ \exists \,\,J_k
\text{ as in (c)}}}  \! |I_s|\\
\le &\,\, C_n\, |F|^{\f{2}{q}} \sum_{m \ge 1} \sum_{\substack{s\in T\\ 2^{-mn}|J_{k(s)}|= |I_s| }} \! 
2^{-mn}|J_{k(s)}|\\
\le &\,\, C_n\, |F|^{\f{2}{q}} \sum_{m \ge 1} (2^m\times (\# \textup{ of sides}) + 2^n)2^{-mn}
 \sum_{k}|J_k|\\
\le &\,\, C_n\, |F|^{\f{2}{q}} |I_t|, 
\end{align*}
where we have used the disjointness of the $J_k$'s.
This finishes case (c).

We now consider case (a). Using (\ref{122}) we can write 
$$
\bigg(\sum_{s\in T} \Big|\sum_{k:\,\, J_k\text{ as in (a)}}
\langle b_k, \phi_s\rangle \Big|^2\bigg)^{\!\f12}
\le C_\gamma \, |F|^{\f1q}  \bigg(
\sum_{s\in T} \Big|\sum_{k:\,\, J_k\subseteq 3I_s }
|J_k |^{\f32} \f{|J_k|^{\f12}}{|I_s|^{\f32}}   \Big|^2  \bigg)^{\!\f12}
$$
and we control the expression inside the parenthesis above by 
$$
   \sum_{s\in T}  \bigg(\sum_{k:\,\, J_k\subseteq 3I_s }
|J_k |^3 \bigg)\bigg(\sum_{k:\,\, J_k\subseteq 3I_s }
\f{|J_k|}{|I_s|^{3}} \bigg)  
\le   \sum_{k:\,\, J_k\subseteq 3I_t} |J_k|^{3}
\sum_{\substack{s\in T\\ J_k\subseteq 3I_s}}   
\f{1}{|I_s|^{2}}
$$
in view of the Cauchy-Schwarz inequality and of the fact that the 
dyadic cubes $J_k$ are disjoint and contained in $3I_s$. 
Finally note that the last sum 
above adds up to at most $C_n\, |J_k|^{-2}$ since for every 
dyadic cube $J_k$ there exist at most $2^n+1+(\# \textup{ of sides})$
 dyadic cubes of 
a given size whose triples contain it.  The required estimate 
$C_{n,\gamma}\, |F|^{\f1q} \, |I_t|^{\f12}$ now follows.

Finally we deal with case  (b) which is the most difficult case. 
We split the set of $k$ into two subsets, those for which 
$J_k\subseteq 3I_t$ and those for which $J_k\nsubseteq 3I_t$, 
(recall $|J_k|\le |I_t|$.) 
Whenever  $J_k\nsubseteq 3I_t$ we have $d(k,s)\approx d(k,t)$. 
In this case we use 
 Minkowski's inequality below and estimate (\ref{132}) with $\gamma >n$ to obtain 
the estimate 
\begin{align*}
\bigg( \sum_{s\in T} \Big|\sum_{k:\,\,J_k\nsubseteq 3I_t} \langle
b_k,\phi_s\rangle
\Big|^2\bigg)^{\!\f12}
\,\, &\le \,\,\sum_{k:\,\,J_k\nsubseteq 3I_t} \bigg( \sum_{s\in T} | \langle
b_k,\phi_s\rangle |^2\bigg)^{\!\f12} \\
\,\, &\le \,\,C_\gamma |F|^{\f1q} \sum_{k:\,\,J_k\nsubseteq 3I_t}|J_k| \bigg( \sum_{s\in
T}   \f{|I_s|^{\f{2\gamma}{n}-1} }{d(k,s)^{2\gamma}} \bigg)^{\!\f12} \\
\,\, &\le \,\,C_\gamma |F|^{\f1q} \sum_{k:\,\,J_k\nsubseteq 3I_t}\f{|J_k|}{d(k,t)^{\gamma}}
\bigg( \sum_{s\in T}   |I_s|^{\f{2\gamma}{n}-1}  \bigg)^{\!\f12} \\
\,\, &\le \,\,C_\gamma |F|^{\f1q}\,|I_t|^{\f{\gamma}{n}-\f12} \sum_{k:\,\,J_k\nsubseteq
3I_t}\f{|J_k|   }{d(k,t)^{\gamma}} \\
\,\, &\le \,\,C_\gamma |F|^{\f1q}\,|I_t|^{\f{\gamma}{n}-\f12} \sum_{l=1}^\nf
\sum_{ k: d(k,t)\approx 2^l|I_t|^{\f1n}}\f{|J_k|  }{(2^l|I_t|^{\f1n})^{\gamma}}  \, . 
\end{align*}
But note that all the $J_k$ with $d(k,t)\approx 2^l|I_t|^{\f1n}$ are contained in 
$2^{l+2}I_t$ and since they are disjoint we can estimate the last 
sum above by $C 2^{lm} |I_t| (2^l |I_t|^{\f1n})^{-\gamma}$. The required 
estimate $C_\gamma |F|^{\f1q} \, |I_t|^{\f12}$ now  follows. 

Next we consider the sum below in which we use estimate (\ref{122})
{\allowdisplaybreaks
\begin{align}
&\,\,\bigg( \sum_{s\in T} \Big|\sum_{\substack{k:\,\,J_k\subseteq 3I_t\\
J_k\cap 3I_s=\emptyset \\ |J_k|\le |I_s|}} \langle
b_k,\phi_s\rangle
\Big|^2\bigg)^{\!\f12} \notag\\
  \le &\,\, C_\gamma\, |F|^{\f1q}
\bigg( \sum_{s\in T} \bigg|\sum_{\substack{k:\,\,J_k\subseteq 3I_t\\
J_k\cap 3I_s=\emptyset \\ |J_k|\le |I_s|}} |J_k|^2 |I_s|^{-\f32 }
\bigg(\f{|I_s|^{\f1n}}{d(k,s)}\bigg)^{\!\gamma}
\bigg|^2\bigg)^{\!\f12} \notag \\
\le &\,\, C_\gamma\, |F|^{\f1q}
\bigg\{ \sum_{s\in T} \bigg[\sum_{\substack{k:\,\,J_k\subseteq 3I_t\\
J_k\cap 3I_s=\emptyset \\ |J_k|\le |I_s|}} \f{|J_k|^3 }{|I_s|^{2}}
 \bigg(\f{|I_s|^{\f1n}}{d(k,s)} \bigg)^{ \!\gamma}\bigg]
\bigg[\sum_{\substack{k:\,\,J_k\subseteq 3I_t\\
J_k\cap 3I_s=\emptyset \\ |J_k|\le |I_s|}}  
\f{|J_k|}{|I_s|} \bigg(\f{d(k,s)}{|I_s|^{\f1n}} \bigg)^{\! -\gamma} 
\bigg]\bigg\}^{\!\f12} \lab{997}. 
\end{align}}
The second sum above can be estimated by
$$
 \sum_{\substack{k:\,\,J_k\subseteq 3I_t\\
J_k\cap 3I_s=\emptyset \\ |J_k|\le |I_s|}}  \int_{J_k} 
 \bigg(\f{|x-c(I_s)|}{|I_s|^{\f1n}} \bigg)^{\! -\gamma}   \f{dx}{|I_s|}
 \le \int_{(3I_s)^c} 
 \bigg(\f{|x-c(I_s)|}{|I_s|^{\f1n}} \bigg)^{\! -\gamma}  \f{ dx}{|I_s|}
\le C_\gamma.
$$
Putting this estimate into (\ref{997}), we have
{\allowdisplaybreaks
\begin{align*}
\le &\,\, C_\gamma\, |F|^{\f1q}
\bigg\{ \sum_{s\in T} \sum_{\substack{k:\,\,J_k\subseteq 3I_t\\
J_k\cap 3I_s=\emptyset \\ |J_k|\le |I_s|}} |J_k|^3 |I_s|^{-2} 
\bigg(\f{|I_s|^{\f1n}}{d(k,s)} \bigg)^{\! \gamma}
\bigg\}^{\!\f12} \\
\le &\,\, C_\gamma\, |F|^{\f1q}
\bigg\{  \sum_{\substack{k:\,\,J_k\subseteq 3I_t\\
J_k\cap 3I_s=\emptyset \\ |J_k|\le |I_s|}} |J_k|^3
\sum_{m \ge \f{\log |J_k|}{n}}2^{-2mn} \sum_{\substack{s \in T \\ |I_s|=2^{mn}}}
 \bigg(\f{d(k,s)}{2^m} \bigg)^{ -\gamma}
\bigg\}^{\!\f12} \\
\le &\,\, C_\gamma\, |F|^{\f1q}
\bigg\{  \sum_{\substack{k:\,\,J_k\subseteq 3I_t\\
J_k\cap 3I_s=\emptyset \\ |J_k|\le |I_s|}} |J_k|^3
\sum_{m \ge \f{\log |J_k|}{n}}2^{-2mn} \bigg\}^{\!\f12} \\
\le &\,\, C_\gamma\, |F|^{\f1q}
\bigg\{  \sum_{\substack{k:\,\,J_k\subseteq 3I_t\\
J_k\cap 3I_s=\emptyset \\ |J_k|\le |I_s|}} |J_k|^3
|J_k|^{-2} \bigg\}^{\!\f12} \\
\le &\,\, C_\gamma\, |F|^{\f1q}\,
|I_t|^{\!\f12}  \, . 
\end{align*}
} $\!\!\!$

There is also the subcase of case (b) in which $|J_k|\ge |I_s|$. 
Here we have the two special subcases: $I_s\cap 3J_k=\emptyset$
and $I_s\subseteq 3J_k=\emptyset$. We begin with the first
of these special subcases in which we use estimate
(\ref{142}). We have 
{\allowdisplaybreaks
\begin{align}
&\,\,\bigg( \sum_{s\in T} \Big|\sum_{
\substack{k:\,\,J_k\subseteq 3I_t\\
J_k\cap 3I_s=\emptyset \\ |J_k|> |I_s|\\
I_s\cap 3J_k=\emptyset}} \langle
b_k,\phi_s\rangle
\Big|^2\bigg)^{\!\f12}\lab{995}\\
  \le &\,\, C_\gamma\, |F|^{\f1q}
\bigg( \sum_{s\in T} \Big|\sum_{\substack{k:\,\,J_k\subseteq 3I_t\\
J_k\cap 3I_s=\emptyset \\ |J_k|> |I_s|\\
I_s\cap 3J_k=\emptyset }}  |I_s|^{ \f12 }
\f{|I_s|^{\f{\gamma}{n}}}{d(k,s)^\gamma}\Big|^2 \bigg)^{\!\f12} \notag\\
\le &\,\, 
C_\gamma\, |F|^{\f1q}
\bigg( \sum_{s\in T} \bigg[\sum_{\substack{k:\,\,J_k\subseteq 3I_t\\
J_k\cap 3I_s=\emptyset \\ |J_k|> |I_s|\\
I_s\cap 3J_k=\emptyset }}  \f{|I_s|^{ 2}}{|J_k|}
\f{|I_s|^{\f{\gamma}{n}}}{d(k,s)^\gamma}\bigg]
\bigg[\sum_{\substack{k:\,\,J_k\subseteq 3I_t\\
J_k\cap 3I_s=\emptyset \\ |J_k|> |I_s|\\
I_s\cap 3J_k=\emptyset }} \f{|J_k|}{|I_s| }
\f{|I_s|^{\f{\gamma}{n}}}{d(k,s)^\gamma}\bigg]
 \bigg)^{\!\f12}\lab{996}.
\end{align}}
 Since $I_s \cap 3J_k = \emptyset$ we have that 
$d(k,s)\approx |x-c(I_s)|$ for every $x\in J_k$. Therefore 
the second term inside square brackets above satisfies
\begin{equation*}
\sum_{\substack{k:\,\,J_k\subseteq 3I_t\\
J_k\cap 3I_s=\emptyset \\ |J_k|> |I_s|\\
I_s\cap 3J_k=\emptyset }} \f{|J_k|}{|I_s| }
\f{|I_s|^{\f{\gamma}{n}}}{d(k,s)^\gamma}\le 
\sum_{k}  
\int_{J_k} 
 \Big(\f{|x-c(I_s)|}{|I_s|^{\f1n}} \Big)^{\! -\gamma}   \f{dx}{|I_s|} \le C_\gamma.
\end{equation*}
Putting this estimate into (\ref{996}), we obtain
{\allowdisplaybreaks
\begin{align*}
&\,\, C_\gamma |F|^{\f1q} \bigg(\sum_{s\in T} \sum_{\substack{k:\,\,J_k\subseteq 3I_t\\
J_k\cap 3I_s=\emptyset \\ |J_k|> |I_s|\\
I_s\cap 3J_k=\emptyset }}  \f{|I_s|^{ 2}}{|J_k|}
\f{|I_s|^{\f{\gamma}{n}}}{d(k,s)^\gamma} \bigg)^{\!\f12}\\
\le &\,\, C_\gamma   |F|^{\f1q}\bigg(\sum_{s\in T}
|I_s|\sum_{\substack{k:\,\,J_k\subseteq 3I_t\\ J_k\cap 3I_s=\emptyset \\ |J_k|> |I_s|\\
I_s\subseteq 3J_k }} 
\f{|I_s|^{\f{\gamma}{n}}}{d(k,s)^\gamma} \bigg)^{\!\f12}\\
\le &\,\, C_\gamma   |F|^{\f1q}\bigg(
\sum_{k: J_k\subseteq 3I_t} |J_k|
\sum_{m=0}^\nf 2^{-mn}\!\!\!\!\!\!\!
\sum_{\substack{s:\, I_s\subseteq
3J_k\\ J_k\cap 3I_s=\emptyset \\
|I_s|=2^{-mn}|J_k|}}\!\!\!\f{|I_s|^{\f{\gamma}{n}}}{d(k,s)^\gamma} \bigg)^{\!\f12}.
\end{align*}}
Since the last sum above is at most a constant  (\ref{995})  satisfies
the estimate $C_\gamma \,|F|^{\f1q}\, |I_t|^{\f12}$.

Finally there is  the subcase of case (b) in which $|J_k|\ge
|I_s|$  
and $I_s\subseteq 3J_k=\emptyset$. Here again we use estimate
(\ref{142}).  We have 
\begin{equation}\lab{0995}
 \bigg( \sum_{s\in T} \Big|\sum_{
\substack{k:\,\,J_k\subseteq 3I_t\\
J_k\cap 3I_s=\emptyset \\ |J_k|> |I_s|\\
I_s\subseteq 3J_k }} \langle
b_k,\phi_s\rangle
\Big|^2\bigg)^{\!\f12} 
  \le   C_\gamma\, |F|^{\f1q}
\bigg( \sum_{s\in T}  |I_s|  \,\Big|\sum_{\substack{k:\,\,J_k\subseteq
3I_t\\ J_k\cap 3I_s=\emptyset \\ |J_k|> |I_s|\\
I_s\subseteq 3J_k }} 
\f{|I_s|^{\f{\gamma}{n}}}{d(k,s)^\gamma}\Big|^2 \bigg)^{\!\f12}  \, . 
\end{equation}
Let us  make some   observations.  For a fixed $s$
there exists at most finitely many $J_k$'s contained in 
$3I_t$ with size at least $|I_s|$. 
Consider the following sets for $\alpha\in\{0,1,2, \ldots\}$,
$$
\mathcal J^\alpha:=\{J_k \textup{ as in the sum above}: 
2^\alpha|I_s|^{\f1n} \le d(k,s) < 2^{\alpha+1}|I_s|^{\f1n}\}.
$$ 
We would like to know that for all $\alpha$ 
the cardinality of $\mathcal J^\alpha$ is bounded
by a fixed constant depending only on dimension.  
This would allow us to work with a single cube
$J^\alpha(s)$ from each set at the cost of a constant in the sum below.
Fix $\alpha\in\{0,1,2, \ldots\}$ and note that $I_s\subseteq 3J_k$ and $d(k,s)
> 2^\alpha|I_s|^{\f1n}$ implies that
$|J_k| > 2^{\alpha n}|I_s|$.  It is clear that
the cardinality of $\mathcal J^{\alpha}$ would be largest if we had 
$|J_k| = 2^{\alpha+1}|I_s|$ for all $J_k \in \mathcal J^\alpha$.  Then  the
cube of size $7^n 2^{\alpha n} |I_s|$ centered at $I_s$ would contain all elements of
$\mathcal J_k$.  This bounds the number of such elements by
 $\big(\f{7}{2}\big)^n$.

Then using the Cauchy-Schwarz inequality we obtain 
\begin{align*}
\Big|\sum_{\substack{k:\,\,J_k\subseteq
3I_t\\ J_k\cap 3I_s=\emptyset \\ |J_k|> |I_s|\\
I_s\subseteq 3J_k }} 
\f{|I_s|^{\f{\gamma}{n}}}{d(k,s)^{\gamma}}\Big|^2  \,\, 
\le &\,\, 
\Big(\f{7}{2}\Big)^{2n}\Big|\sum_{\alpha=1}^\nf  
\f{|I_s|^{\f{\gamma}{2n}}}{\textup{dist}(J^{\alpha}(s), I_s)^{\f
{\gamma}{2}}}\f{1}{2^{\f{\alpha\gamma}{2}}}
\Big|^2  \\
\,\, \le &\,\, C_n
 \sum_{\alpha=1}^\nf  
\f{|I_s|^{ \f{\gamma}{n} }}{\textup{dist}(J^\alpha(s), I_s)^{\gamma} }  \\
\,\, \le &\,\, C_n 
\sum_{\substack{k:\,\,J_k\subseteq
3I_t\\ J_k\cap 3I_s=\emptyset \\ |J_k|> |I_s|\\
I_s\subseteq 3J_k }} 
\f{|I_s|^{\f{\gamma}{n}}}{d(k,s)^\gamma} 
\end{align*}
Putting  this estimate into the right hand side of  (\ref{0995}), 
  the estimate $C_{n,\gamma} \,|F|^{\f1q}\, |I_t|^{\f12}$ now follows as
 in the previous case.
This concludes the proof of Lemma \ref{L}. 

\section{applications}

We conclude by discussing some applications. We show how one 
can strengthen the results of the previous sections to obtain  
distributional estimates for the function $\mathcal D_r(\chi_F)$
similar to those in the paper of  Sj\"olin \cite{sjolin}. 

We showed in section 4 that for any measurable set $E$ there is 
a set $E'$ of at least half the measure of $E$ such that 
\begin{equation}\lab{sjolin}
\Big| \int_{E'} \mathcal D_r(\chi_F)\, dx \Big| \le C\, 
\min(|E|,|F|)\Big(  1+\Big|\log \f{|F|}{|E|} \Big|\Big)
\end{equation}
for some constant $C$  depending only on the dimension. 
For $\la>0$   we define  
$$
E_\la=\big\{|\mathcal D_r(\chi_F)|>\la\big\}
$$
and also 
\begin{align*}
E_\la^1=&\big\{ \text{Re }\mathcal D_r(\chi_F) >\la\big\}   
&&E_\la^2=\big\{ \text{Re }\mathcal D_r(\chi_F) <-\la\big\}\\
E_\la^3=& \big\{ \text{Im }\mathcal D_r(\chi_F) > \la\big\} 
&&E_\la^4=\big\{ \text{Im }\mathcal D_r(\chi_F) <-\la\big\}. 
\end{align*}
We apply (\ref{sjolin}) to each set $E_\la^j$ to obtain 
$$
\la \, |E_\la^j| \le \min(|E_\la^j|,|F|)\Big(  1+\Big|\log \f{|F|}{|E_\la^j|} \Big|\Big).
$$
Using this fact in combination with the easy observation that for $a>1$
$$
\f{a}{\log a}\le \f{1}{\la} \implies a \le \f{10 }{\la}\log( \f{1}{\la})\, , 
$$
to obtain that 
$$
|E_\la^j| \le C'|F|\begin{cases}
 \f{1}{\la}\log( \f{1}{\la})   &\text{when $\la<\f12$}\\
e^{-c\la}  &\text{when $\la\ge \f12$.} 
 \end{cases}
$$
Since $|E_{2\sqrt{2}\,\la}|\le \sum_{j=1}^4 |E_\la^j|$ we conclude 
a similar estimate for $E_\la$.

Next we obtain  similar distributional estimates for maximally
modulated singular integrals $\mathcal M$ such as the 
maximally modulated Hilbert transform (i.e.
Carleson's operator) or the
maximally modulated Riesz transforms
$$
\mathcal M (f)(x)=\sup_{\xi\in\rn}\bigg|
\int_{\rn} \f{x_j-y_j}{|x-y|^{n+1}} e^{2\pi i \xi\cdot y}f(y)\, dy\, \bigg|\, . 
$$ 
To achieve this in the one dimensional setting, one
 applies an averaging argument similar to that in \cite{LT}
 to both terms of estimate (\ref{sjolin}) to recover a similar estimate  
with the Carleson operator.  For more general homogeneous 
singular integrals 
 with sufficiently smooth kernels, one applies the averaging argument
 to  suitable modifications 
 of the operators $\mathcal D_r$ as in \cite{em}.
Then one obtains a version of estimate  
(\ref{sjolin}) in which $\mathcal D_r(\chi_F)$ is replaced by 
 $\mathcal M(\chi_F)$. The same procedure as above then 
yields the distributional estimate
$$
\big|\big\{|\mathcal M(\chi_F)|>\la\big\} \big| \le C_n'|F|\begin{cases}
 \f{1}{\la}\log( \f{1}{\la})    &\text{when $\la<\f12$}\\
e^{-c\la}  &\text{when $\la\ge \f12$.} 
 \end{cases}
$$
which recovers Lemma 1.2 in \cite{sjolin}. It should be noted
that the corresponding estimate 
\begin{equation}\label{last-sjolin}
\big|\big\{|\mathcal D_r(\chi_F)|>\la\big\} \big| \le C_n|F|\begin{cases}
 \f{1}{\la}\log( \f{1}{\la})    &\text{when $\la<\f12$}\\
e^{-c\la}  &\text{when $\la\ge \f12$.} 
 \end{cases}
\end{equation}
 obtained here for $\mathcal D_r$ is  
stronger as it concerns an ``unaveraged version'' 
of all the aforementioned maximally 
modulated singular integrals $\mathcal M$.

Using the idea employed in Sj\"olin \cite{sjolin0} we can   obtain the following 
result as a consequence of (\ref{last-sjolin}). Let $B$ be a ball in $\rn$. 

\begin{proposition}
(i) If $\int_B |f(x)|\, \log^+|f(x)|\, \log^+ \log^+|f(x)|\, dx<\nf$, then 
$\mathcal D_r(f)$ is finite a.e. on $B$. \\
(ii) If $\int_B |f(x)|\, (\log^+|f(x)|)^2\,  dx<\nf$, then $\mathcal D_r(f)$
is integrable over $B$. \\
(iii) For all $\la>0$ we have 
$$
\big|\{ x\in \rn:\,\, |\mathcal D_r(f)(x)|>\la\}\big| \le 
C\, e^{-c\la/\|f\|_{\li}}
$$
where $C,c$ only depend on the dimension (in particular they are independent of 
 the measurable function $N:\rn\to \rn$.)
\end{proposition}

\end{document}